\documentclass[runningheads,svgnames]{llncs}
\usepackage[T1]{fontenc}
\usepackage{graphicx}
\usepackage{pgfplots}
\usepgfplotslibrary{statistics}
\usepgfplotslibrary{fillbetween}

\usepackage{amsmath,amsfonts,amssymb}
\usepackage{algpseudocode}
\usepackage{algorithm} 
\usepackage{tikz,tikz-3dplot}
\usepackage{bm}
\usepackage{multirow}

\usepackage[colorlinks=true, allcolors=blue]{hyperref}
\usepackage{bookmark}
\errorcontextlines\maxdimen

\makeatletter
\newcommand*{\algrule}[1][\algorithmicindent]{\makebox[#1][l]{\hspace*{.5em}\thealgruleextra\vrule height \thealgruleheight depth \thealgruledepth}}%
\newcommand*{\thealgruleextra}{}
\newcommand*{\thealgruleheight}{.75\baselineskip}
\newcommand*{\thealgruledepth}{.25\baselineskip}

\newcount\ALG@printindent@tempcnta
\def\ALG@printindent{%
	\ifnum \theALG@nested>0% is there anything to print
	\ifx\ALG@text\ALG@x@notext% is this an end group without any text?
	\else
	\unskip
	\addvspace{-1pt}% FUDGE to make the rules line up
	\ALG@printindent@tempcnta=1
	\loop
	\algrule[\csname ALG@ind@\the\ALG@printindent@tempcnta\endcsname]%
	\advance \ALG@printindent@tempcnta 1
	\ifnum \ALG@printindent@tempcnta<\numexpr\theALG@nested+1\relax% can't do <=, so add one to RHS and use < instead
	\repeat
	\fi
	\fi
}%
\usepackage{etoolbox}
\patchcmd{\ALG@doentity}{\noindent\hskip\ALG@tlm}{\ALG@printindent}{}{\errmessage{failed to patch}}
\makeatother

\newbox\statebox
\newcommand{\myState}[1]{%
	\setbox\statebox=\vbox{#1}%
	\edef\thealgruleheight{\dimexpr \the\ht\statebox+1pt\relax}%
	\edef\thealgruledepth{\dimexpr \the\dp\statebox+1pt\relax}%
	\ifdim\thealgruleheight<.75\baselineskip
	\def\thealgruleheight{\dimexpr .75\baselineskip+1pt\relax}%
	\fi
	\ifdim\thealgruledepth<.25\baselineskip
	\def\thealgruledepth{\dimexpr .25\baselineskip+1pt\relax}%
	\fi
	\State #1%
	\def\thealgruleheight{\dimexpr .75\baselineskip+1pt\relax}%
	\def\thealgruledepth{\dimexpr .25\baselineskip+1pt\relax}%
}
\usepackage{comment}

\usepackage{color}

\begin{document}
\title{An Integer Linear Programming Model\\ for the Evolomino Puzzle}
\author{Andrei V. Nikolaev\orcidID{0000-0003-4705-2409}\\ \and
Yuri A. Myasnikov\orcidID{0009-0007-1571-0071}}
\authorrunning{A.V. Nikolaev and Y.A. Myasnikov}
\institute{P.G. Demidov Yaroslavl State University, Yaroslavl, Russia
	\email{andrei.v.nikolaev@gmail.com, ya\_myasnikov@vk.com}
}
\maketitle              % typeset the header of the contribution
\begin{abstract}
Evolomino is a pencil‑and‑paper logic puzzle published by the Japanese company Nikoli, renowned for culture‑independent puzzles such as Sudoku, Kakuro, and Slitherlink.
Its name reflects the core mechanic: the polyomino-like blocks drawn by the player must gradually ``evolve'' according to the directions indicated by arrows pre-printed on a rectangular grid.

In this paper, we formalize the rules of Evolomino as an integer linear programming (ILP) model, encoding block evolution, connectivity, and consistency requirements through linear constraints.
Furthermore, we introduce an algorithm for generating random Evolomino instances, utilizing this ILP framework to ensure solution uniqueness.
Computational experiments on a custom benchmark dataset demonstrate that a state-of-the-art CP-SAT solver successfully handles puzzle instances of up to $10 \times 10$ within one second and up to $18 \times 18$ within one minute.

%The abstract should briefly summarize the contents of the paper in

\keywords{Evolomino \and Pencil-and-paper logic puzzle \and Polyominoes \and Integer linear programming  \and Puzzle generation \and Solution uniqueness.}
	
%	First keyword  \and Second keyword \and Another keyword.}
\end{abstract}

\section{Introduction}
	
	Nikoli is a renowned Japanese publisher specializing in board games and, in particular, logic puzzles. Founded in 1980 by Maki Kaji, often referred to as ``the father of Sudoku'', the company rose to worldwide prominence with the global spread of Sudoku.
	
	Nikoli is especially known for its extensive library of culture‑independent puzzles that do not rely on any specific language or alphabet and are often purely logical or numerical in nature.
	The company’s namesake magazine was the first puzzle magazine in Japan, and over the years it has introduced numerous new puzzle genres and brought several new games to the Japanese audience.
	
	In this paper, we study the relatively new pencil‑and‑paper logic puzzle Evolomino, which first appeared in \textit{Puzzle Communication Nikoli}, Vol.~182, in March 2023~\cite{PuzzleCommunicationNikoli182}.

	Here we present the rules of the Evolomino puzzle as they appear on Nikoli’s official website~\cite{Evolomino_Nikoli_web_page}:
	
	\begin{itemize}
		\item Evolomino is played on a rectangular board with white and shaded cells. Some white cells contain pre-drawn squares and arrows.
		
		\item The player must draw squares ($\square$) in some of the white cells.
		
		\item A polyomino-like group of squares connected vertically and horizontally is called a \textit{block} (including only one square). Each block must contain exactly one square placed on a pre-drawn arrow.
		
		\item Each arrow must pass through at least two blocks.
		
		\item The second and later blocks on the route of an arrow from start to finish must progress by adding one square to the previous block without rotating or flipping.
	\end{itemize}
	
	An example of the Evolomino puzzle and its solution from Nikoli's official website~\cite{Evolomino_Nikoli_web_page} is shown in Fig.~\ref{Fig_Evolomino_example}.

	\begin{figure}[t]
	\centering
	\begin{tikzpicture}[scale=0.7]
		
		\foreach \x in {0,...,5}
		\draw (\x+0.5,0.5)--(\x+0.5,5.5);
		
		\foreach \y in {0,...,5}
		\draw (0.5,\y+0.5)--(5.5,\y+0.5);

		\foreach \x in {5}	
		\foreach \y in {1}
		\draw [fill = MidnightBlue!30] (\x-0.5,\y-0.5) -- (\x+0.5,\y-0.5) -- (\x+0.5,\y+0.5) -- (\x-0.5,\y+0.5) -- cycle;	
		
		\foreach \x in {4}	
		\foreach \y in {4}
		\draw [fill = MidnightBlue!30] (\x-0.5,\y-0.5) -- (\x+0.5,\y-0.5) -- (\x+0.5,\y+0.5) -- (\x-0.5,\y+0.5) -- cycle;

		\draw [semithick,->,>=stealth] (1,1) -- (3+0.1,1);
		\draw [semithick,->,>=stealth] (1,3) -- (1,5) -- (4+0.1,5);
		\draw [semithick,->,>=stealth] (5,2) -- (4,2) -- (4,3+0.1);
		
		\foreach \x in {1}
		\foreach \y in {3}
		\draw (\x-0.25,\y-0.25) -- (\x+0.25,\y-0.25) -- (\x+0.25,\y+0.25) -- (\x-0.25,\y+0.25) -- cycle;

		\foreach \x in {5}
		\foreach \y in {4}
		\draw (\x-0.25,\y-0.25) -- (\x+0.25,\y-0.25) -- (\x+0.25,\y+0.25) -- (\x-0.25,\y+0.25) -- cycle;
		
		\node at (3,0) {Sample puzzle};
		
		\begin{scope}[xshift=8cm]
			
			\foreach \x in {0,...,5}
			\draw (\x+0.5,0.5)--(\x+0.5,5.5);
			
			\foreach \y in {0,...,5}
			\draw (0.5,\y+0.5)--(5.5,\y+0.5);

			\foreach \x in {5}	
			\foreach \y in {1}
			\draw [fill = MidnightBlue!30] (\x-0.5,\y-0.5) -- (\x+0.5,\y-0.5) -- (\x+0.5,\y+0.5) -- (\x-0.5,\y+0.5) -- cycle;	
			
			\foreach \x in {4}	
			\foreach \y in {4}
			\draw [fill = MidnightBlue!30] (\x-0.5,\y-0.5) -- (\x+0.5,\y-0.5) -- (\x+0.5,\y+0.5) -- (\x-0.5,\y+0.5) -- cycle;

			\draw [semithick,->,>=stealth] (1,1) -- (3+0.1,1);
			\draw [semithick,->,>=stealth] (1,3) -- (1,5) -- (4+0.1,5);
			\draw [semithick,->,>=stealth] (5,2) -- (4,2) -- (4,3+0.1);
			
			\foreach \x in {1}
			\foreach \y in {1,3}
			\draw (\x-0.25,\y-0.25) -- (\x+0.25,\y-0.25) -- (\x+0.25,\y+0.25) -- (\x-0.25,\y+0.25) -- cycle;
			
			\foreach \x in {2}
			\foreach \y in {4,5}
			\draw (\x-0.25,\y-0.25) -- (\x+0.25,\y-0.25) -- (\x+0.25,\y+0.25) -- (\x-0.25,\y+0.25) -- cycle;
			
			\foreach \x in {3}
			\foreach \y in {1,3}
			\draw (\x-0.25,\y-0.25) -- (\x+0.25,\y-0.25) -- (\x+0.25,\y+0.25) -- (\x-0.25,\y+0.25) -- cycle;
			
			\foreach \x in {4}
			\foreach \y in {1,3,5}
			\draw (\x-0.25,\y-0.25) -- (\x+0.25,\y-0.25) -- (\x+0.25,\y+0.25) -- (\x-0.25,\y+0.25) -- cycle;
			
			\foreach \x in {5}
			\foreach \y in {2,4,5}
			\draw (\x-0.25,\y-0.25) -- (\x+0.25,\y-0.25) -- (\x+0.25,\y+0.25) -- (\x-0.25,\y+0.25) -- cycle;
			
			\node at (3,0) {Solution};
		\end{scope}
	\end{tikzpicture}
	\caption {Example of an Evolomino puzzle}
	\label {Fig_Evolomino_example}
\end{figure}

As a relatively new puzzle, Evolomino has not yet been extensively studied from a theoretical computer science perspective.  
In~\cite{Nikolaev2025}, a reduction from 3‑SAT was used to prove that deciding whether an Evolomino puzzle has at least one valid solution is NP‑complete. Furthermore, since the reduction is parsimonious and preserves the number of distinct solutions, it follows that the associated counting problem is $\texttt{\#}$P‑complete.

In this paper, we formulate the rules of the Evolomino puzzle as an integer linear programming (ILP) model. ILP models are among the most established and versatile methods for solving not only combinatorial optimization problems but also a wide range of puzzles and logic challenges.
Examples include MILP formulations for Battleship~\cite{Meuffels2010}, the Eternity puzzle~\cite{Burkardt2023}, Nonograms (both black‑and‑white~\cite{Bosch2001} and colored~\cite{Mingote2009}), Kakuro~\cite{Simonis2008}, Shakashaka~\cite{Demaine2014}, Smashed Sums~\cite{Kececi2021}, the Survo puzzle~\cite{Sungur2025}, Tantrix~\cite{Kino2012}, and even letter‑substitution cipher attacks~\cite{Ravi2009}.
ILP‑based approaches have also been proposed for several well‑known Nikoli puzzles, including Sudoku~\cite{Bartlett2008}, Kakuro~\cite{Simonis2008}, and Futoshiki~\cite{Sungur2022}.
The book by Tony Hürlimann~\cite{Hurlimann2025} offers a comprehensive overview of mathematical modeling, featuring more than 140 puzzles as illustrative case studies.

Although integer programming is one of Karp's original 21 NP-complete problems~\cite{Karp1972}, many algorithms have proven remarkably effective in practice. For example, Koch et al.~\cite{Koch2022} demonstrated that between 2001 and 2020 computer hardware became roughly $20$ times faster and algorithms improved by a factor of about $50$, resulting in an overall speed-up of approximately $1\,000$ for MILP models. Moreover, many instances considered unsolvable twenty years ago can now be solved within seconds.
Therefore, once a puzzle can be formulated as an integer programming problem, we may hope to use modern solvers to handle it efficiently in practice.

The rest of the paper is organized as follows.
Section~2 presents the integer linear programming constraints that model the rules of the Evolomino puzzle.
Section~3 outlines the algorithm for generating random puzzle instances.
Section~4 reports the results of computational experiments conducted on a custom benchmark dataset, produced by the generator from Section~3, with CP-SAT~\cite{cpsatlp} from the OR-Tools~\cite{OR-Tools} library serving as the main solver.

\section{ILP-model}

\subsection{Preliminary notation}

\begin{itemize}
	
	\item $\mathcal{C}$ --- the set of all cells of a rectangular $m \times n$ board. 
	
	For simplicity, the cells of the board are numbered by a single linear index row by row, starting from the top‑left corner.
	
	\item $\mathcal{A}$ --- the set of all arrows. 
	
	\item $\mathcal{P}_a \subset \mathcal{C}$ --- the set of cells along arrow $a \in \mathcal{A}$. 
	
	\item $\operatorname{next}_{a}(i)$ --- the cell that follows $i \in \mathcal{P}_a$ in the direction of arrow $a \in \mathcal{A}$.
	
	\item $\mathcal{P}^{\prec i}_{a}$ --- the set of cells preceding cell $i \in \mathcal{P}_a$ in the direction of arrow $a \in \mathcal{A}$.
	
	\item $K_a$ --- the maximum possible number of blocks along arrow $a \in \mathcal{A}$, computed for each arrow as $K_a = \lceil \frac{|\mathcal{P}_a|}{2} \rceil$.
	
	\item $\mathcal{C}_a$ --- the set of cells that may belong to a polyomino block along arrow $a \in \mathcal{A}$.
	The example of region $\mathcal{C}_a$ is shown in Fig.~\ref{Fig_C_a} (highlighted in green). 
	Starting from an arbitrary cell $i$ on arrow $a$, we perform a breadth-first search to collect all cells $j$ that are reachable from $i$ by a path that does not pass through shaded cells, other arrows, or cells belonging to blocks of other arrows; all such cells are included in $\mathcal{C}_a$.

	\begin{figure}[t]
		\centering
		\begin{tikzpicture}[scale=0.7]
			
			\foreach \x in {0,...,5}
			\draw (\x+0.5,0.5)--(\x+0.5,5.5);
			
			\foreach \y in {0,...,5}
			\draw (0.5,\y+0.5)--(5.5,\y+0.5);
			
	%		\foreach \x in {1,...,5}
	%		\node at (\x,0) {\x};
			
	%		\foreach \y in {1,...,5}
	%		\node at (0,\y) {\y};
			
			\foreach \x in {5}	
			\foreach \y in {1}
			\draw [fill = MidnightBlue!30] (\x-0.5,\y-0.5) -- (\x+0.5,\y-0.5) -- (\x+0.5,\y+0.5) -- (\x-0.5,\y+0.5) -- cycle;	
			
			\foreach \x in {4}	
			\foreach \y in {4}
			\draw [fill = MidnightBlue!30] (\x-0.5,\y-0.5) -- (\x+0.5,\y-0.5) -- (\x+0.5,\y+0.5) -- (\x-0.5,\y+0.5) -- cycle;

			\foreach \x in {1}	
			\foreach \y in {1}
			\draw [fill = green!10] (\x-0.5,\y-0.5) -- (\x+0.5,\y-0.5) -- (\x+0.5,\y+0.5) -- (\x-0.5,\y+0.5) -- cycle;
			
			\foreach \x in {2}	
			\foreach \y in {1,2,4}
			\draw [fill = green!10] (\x-0.5,\y-0.5) -- (\x+0.5,\y-0.5) -- (\x+0.5,\y+0.5) -- (\x-0.5,\y+0.5) -- cycle;
			
			\foreach \x in {3}	
			\foreach \y in {1,2,3,4}
			\draw [fill = green!10] (\x-0.5,\y-0.5) -- (\x+0.5,\y-0.5) -- (\x+0.5,\y+0.5) -- (\x-0.5,\y+0.5) -- cycle;
			
			\foreach \x in {4}	
			\foreach \y in {1}
			\draw [fill = green!10] (\x-0.5,\y-0.5) -- (\x+0.5,\y-0.5) -- (\x+0.5,\y+0.5) -- (\x-0.5,\y+0.5) -- cycle;
			
			\draw [semithick,->,>=stealth] (1,1) -- (3+0.1,1);
			\draw [semithick,->,>=stealth] (1,3) -- (1,5) -- (4+0.1,5);
			\draw [semithick,->,>=stealth] (5,2) -- (4,2) -- (4,3+0.1);
			
			\node at (2,1.25) {\small $a$};
			
			\foreach \x in {1}
			\foreach \y in {3}
			\draw (\x-0.25,\y-0.25) -- (\x+0.25,\y-0.25) -- (\x+0.25,\y+0.25) -- (\x-0.25,\y+0.25) -- cycle;

			\foreach \x in {5}
			\foreach \y in {4}
			\draw (\x-0.25,\y-0.25) -- (\x+0.25,\y-0.25) -- (\x+0.25,\y+0.25) -- (\x-0.25,\y+0.25) -- cycle;
			
		\end{tikzpicture}
		\caption{Example of region $\mathcal{C}_a$ for arrow $a$, highlighted in green}
		\label {Fig_C_a}
	\end{figure}
	
	\item $\operatorname{row}(i)$ and $\operatorname{col}(i)$ --- functions that return the row and column indices of cell $i$, respectively.
	
	\item $M$ --- a sufficiently large constant that exceeds the size of any polyomino block, for example $M = |\mathcal{C}|$.
	
	\item $\mathcal{T}$ --- the set of all possible two-dimensional translations (shifts) applicable to two consecutive blocks, where the second block contains a copy of a first block as a sub-shape displaced by $p$ rows and $q$ columns. It consists of all integer pairs $(p,q)$ with $p \in \{-(m-1), \ldots, m-1\}$ and $q \in \{-(n-1), \ldots, n-1\}$ where $|p|+|q| > 1$ (no zero translations or translations to adjacent cells).
	
	\item $\operatorname{shift}(i,p,q) = i + p \cdot n + q$ --- the function that maps linear index $i$ to the index of the cell shifted by $p$ rows and $q$ columns.
	
	\item $\mathcal{T}^{a}_i \subseteq \mathcal{T}$ --- the set of admissible translations $(p,q)$ from a cell $i \in \mathcal{C}_a$ along arrow $a \in \mathcal{A}$ such that the shift remains within the board:
	$1 \leq \operatorname{row}(i) + p \leq m$, 
	$1 \leq \operatorname{col}(i) + q \leq n$, 
	and $\operatorname{shift}(i,p,q) \in \mathcal{C}_a$ (valid translation inside region $\mathcal{C}_a$).
	An example of the set $\mathcal{T}^{a}_i$ is shown in Fig.~\ref{Fig_T^a_i}.

	\begin{figure}[t]
		\centering
		\begin{tikzpicture}[scale=0.7]
			
			\foreach \x in {0,...,5}
			\draw (\x+0.5,0.5)--(\x+0.5,5.5);
			
			\foreach \y in {0,...,5}
			\draw (0.5,\y+0.5)--(5.5,\y+0.5);
			
			%		\foreach \x in {1,...,5}
			%		\node at (\x,0) {\x};
			
			%		\foreach \y in {1,...,5}
			%		\node at (0,\y) {\y};
			
			\foreach \x in {5}	
			\foreach \y in {1}
			\draw [fill = MidnightBlue!30] (\x-0.5,\y-0.5) -- (\x+0.5,\y-0.5) -- (\x+0.5,\y+0.5) -- (\x-0.5,\y+0.5) -- cycle;	
			
			\foreach \x in {4}	
			\foreach \y in {4}
			\draw [fill = MidnightBlue!30] (\x-0.5,\y-0.5) -- (\x+0.5,\y-0.5) -- (\x+0.5,\y+0.5) -- (\x-0.5,\y+0.5) -- cycle;			
			
			\foreach \x in {1}	
			\foreach \y in {1}
			\draw [fill = orange!15] (\x-0.5,\y-0.5) -- (\x+0.5,\y-0.5) -- (\x+0.5,\y+0.5) -- (\x-0.5,\y+0.5) -- cycle;

	%		\foreach \x in {1}	
	%		\foreach \y in {1}
	%		\draw [fill = green!10] (\x-0.5,\y-0.5) -- (\x+0.5,\y-0.5) -- (\x+0.5,\y+0.5) -- (\x-0.5,\y+0.5) -- cycle;
			
			\foreach \x in {2}	
			\foreach \y in {2,4}
			\draw [fill = green!10] (\x-0.5,\y-0.5) -- (\x+0.5,\y-0.5) -- (\x+0.5,\y+0.5) -- (\x-0.5,\y+0.5) -- cycle;
			
			\foreach \x in {3}	
			\foreach \y in {1,2,3,4}
			\draw [fill = green!10] (\x-0.5,\y-0.5) -- (\x+0.5,\y-0.5) -- (\x+0.5,\y+0.5) -- (\x-0.5,\y+0.5) -- cycle;
			
			\foreach \x in {4}	
			\foreach \y in {1}
			\draw [fill = green!10] (\x-0.5,\y-0.5) -- (\x+0.5,\y-0.5) -- (\x+0.5,\y+0.5) -- (\x-0.5,\y+0.5) -- cycle;
			
			\draw [semithick,->,>=stealth] (1,1) -- (3+0.1,1);
			\draw [semithick,->,>=stealth] (1,3) -- (1,5) -- (4+0.1,5);
			\draw [semithick,->,>=stealth] (5,2) -- (4,2) -- (4,3+0.1);
			
			\node at (1,1.25) {\small $i$};
			\node at (2,1.25) {\small $a$};
			
			\foreach \x in {1}
			\foreach \y in {3}
			\draw (\x-0.25,\y-0.25) -- (\x+0.25,\y-0.25) -- (\x+0.25,\y+0.25) -- (\x-0.25,\y+0.25) -- cycle;

			\foreach \x in {5}
			\foreach \y in {4}
			\draw (\x-0.25,\y-0.25) -- (\x+0.25,\y-0.25) -- (\x+0.25,\y+0.25) -- (\x-0.25,\y+0.25) -- cycle;
			
			\node at (2,2) {\scalebox {0.625} {$(-1,1)$}};
			\node at (2,4) {\scalebox {0.625} {$(-3,1)$}};
			\node at (3,1.25) {\scalebox {0.625} {$(0,2)$}};
			\node at (3,2) {\scalebox {0.625} {$(-1,2)$}};
			\node at (3,3) {\scalebox {0.625} {$(-2,2)$}};
			\node at (3,4) {\scalebox {0.625} {$(-3,2)$}};
			\node at (4,1) {\scalebox {0.625} {$(0,3)$}};
			
		\end{tikzpicture}
		\caption{Example of set $\mathcal{T}^a_i$ of admissible translations $(p,q)$ from a cell $i$ highlighted in orange along arrow $a$; admissible target cells are highlighted in green}
		\label {Fig_T^a_i}
	\end{figure}
	
	\item $\mathcal{T}_{a} = \bigcup_{i \in \mathcal{P}_a} \mathcal{T}^{a}_i$ ---  the set of all admissible two-dimensional translations $(p,q)$ for the blocks associated with arrow $a \in \mathcal{A}$.
	
\end{itemize}

\subsection{Basic variables}

\begin{itemize}
	\item For each cell $i \in \mathcal{C}$ we introduce the \textit{cell-variable}:
	\[
	x_{i} = \begin{cases}
		1,& \text{if cell } i \text{ contains a square } \square,\\
		0,& \text{otherwise}.
	\end{cases}
	\]

	The variables $x_i$ are used to encode the initial puzzle grid: pre‑placed squares ($x_i = 1$) and shaded cells ($x_i = 0$), and they will also store the final solution of the puzzle.
	
%	For simplicity, throughout the text we use a single index for each cell.
	
	\item For each arrow $a \in \mathcal{A}$, each block $k \leq K_a$, and each cell $i \in \mathcal{C}_a$ we introduce the \textit{block-variable}:
	\[
	y_{i}^{ak} = \begin{cases}
		1,& \text{if cell } i \text{ belongs to block } k \text{ on arrow } a,\\
		0,& \text{otherwise}.
	\end{cases}
	\]

	Using the variables $y_{i}^{ak}$ we encode all blocks by associating them with specific cells of the board.
	
	\item For each arrow $a \in \mathcal{A}$ and each block $k \leq K_a$ we introduce the \textit{block activation variable}:
	\[
	b_{ak} = \begin{cases}
		1,& \text{if block } k \text{ on arrow } a \text{ exists},\\
		0,& \text{otherwise}.
	\end{cases}
	\]
	
	Using the variables $b_{ak}$ we encode the number of blocks along each arrow. Since any arrow passes through at least two blocks, we set $b_{a1} = b_{a2} = 1$ for every arrow $a \in \mathcal{A}$.
	
\end{itemize}

\subsection{Basic constraints}

\begin{itemize}

	\item Each square must be assigned to exactly one block:
	\begin{align}
		\sum_{\substack{a \in \mathcal{A} \\ i \in \mathcal{C}_a}}\sum_{k = 1}^{K_a} y_{i}^{ak} &= x_{i} &&\forall i \in \mathcal{C}. \label{constr:one-block-per-cell}
	\end{align}
	If a $\square$ is placed in cell $i$ ($x_i = 1$), then it belongs to exactly one block.

	\item No two consecutive squares may be placed along an arrow:
	\begin{align}
		x_i + x_{\operatorname{next}_a(i)} &\leq 1 && \forall a \in \mathcal{A},\ \forall i \in \mathcal{P}_{a}. \label{constr:no-adjacent-x-along-arrow}
	\end{align}
	
	\item If a block is not activated, then it must be empty:
	\begin{align}
		\sum_{i \in \mathcal{C}_a} y_{i}^{ak} &\leq M \cdot b_{ak}
		&& \forall a \in \mathcal{A},\ k = 1,\ldots, K_a. \label{constr:block-activation}		
	\end{align}
%	where $M$ is a sufficiently large constant, for example the board size $|\mathcal{C}|$.
	
	\item If a block is activated, then it contains exactly one square along the arrow:
	\begin{align}
		\sum_{i \in \mathcal{P}_a} y_{i}^{ak} &= b_{ak}
		&& \forall a \in \mathcal{A},\ k = 1,\ldots, K_a. \label{constr:one-square-per-active-block-on-arrow}
	\end{align}
	
	\item A subsequent block cannot be activated before the previous one:
	\begin{align}
		b_{ak} &\leq b_{a,k-1}
		&& \forall a \in \mathcal{A},\ k = 2,\ldots, K_a. \label{constr:block-activation-order}
	\end{align}
	
	\item The order of blocks must follow the direction of the arrow:
	\begin{align}
		\sum_{j \in \mathcal{P}^{\prec i}_a} y_{j}^{ak}
		&\leq 1 - y^{a,k-1}_{i}
		&& \forall a \in \mathcal{A},\ \forall i \in \mathcal{P}_a,\ k = 2,\ldots, K_a. \label{constr:block-order-along-arrow}
	\end{align}
	Block $k$ cannot appear earlier than block $k-1$ along arrow $a$. If cell $i$ contains the square of block $k-1$ (i.e., $y^{a,k-1}_{i} = 1$), then block $k$ (and, by induction, all subsequent blocks) cannot be placed in any position from $\mathcal{P}^{\prec i}_a$.
	
	\item Two adjacent cells cannot belong to different blocks:
	\begin{align}
		y_{i}^{ak} + y_{i+1}^{a'k'} &\leq 1
		&& \forall a,a' \in \mathcal{A},\ \forall i \in \mathcal{C}_a,\ i+1 \in \mathcal{C}_{a'},\ k,k' = 1,\ldots, K_a, \label{constr:no-adjacent-blocks-horizontal} \\
		y_{i}^{ak} + y_{i+n}^{a'k'} &\leq 1
		&& \forall a,a' \in \mathcal{A},\ \forall i \in \mathcal{C}_a,\ i+n \in \mathcal{C}_{a'},\ k,k' = 1,\ldots, K_a, \label{constr:no-adjacent-blocks-vertical}
	\end{align}
	where $a \neq a'$ or $k \neq k'$, $\operatorname{row}(i) < m$, and $\operatorname{col}(i) < n$.
	
	This restriction is imposed for every pair of blocks and every pair of adjacent cells on the board, excluding cells in the last row and the last column.

	\item For convenience, we introduce the \textit{block size variable} $N_{ak}$ for block $k$ on arrow $a$:
	\begin{align}
		N_{ak} &= \sum_{i \in \mathcal{C}_a} y_{i}^{ak}
		&& \forall a \in \mathcal{A},\ k = 1,\ldots, K_a. \label{constr:N-definition}
	\end{align}
	This variable is not strictly necessary, but it makes the subsequent constraints slightly easier to write.
	
	\item The size of each subsequent block along arrow $a$ must exceed the previous one by exactly one square:
	\begin{align}
		N_{ak} &\geq N_{a,k-1} + 1 - M(1 - b_{ak})
		&& \forall a \in \mathcal{A},\ k = 2,\ldots, K_a, \label{constr:block-size-increment-geq}\\
		N_{ak} &\leq N_{a,k-1} + 1 + M(1 - b_{ak})
		&& \forall a \in \mathcal{A},\ k = 2,\ldots, K_a. \label{constr:block-size-increment-leq}
	\end{align}
	Here, if block $k$ is activated ($b_{ak} = 1$), then constraints
	\eqref{constr:block-size-increment-geq}--\eqref{constr:block-size-increment-leq} reduce to
	\[ N_{ak} = N_{a,k-1} + 1. \]
	Otherwise, if the block is not activated, then
	\eqref{constr:block-size-increment-geq}--\eqref{constr:block-size-increment-leq} reduce to trivial bounds.
	
\end{itemize}

\subsection{Connectivity variables and constraints}

The following group of variables and constraints ensures that each block forms a connected polyomino. To achieve this, we adapt a classic flow-based model by Gavish and Graves~\cite{Gavish1978}.

For each block $k$ on each arrow $a \in \mathcal{A}$ and for every pair of adjacent cells $i$ and $j$ from $\mathcal{C}_a$ (each cell $i$ has up to four neighboring cells, which we denote by $\mathcal{N}(i)$), we introduce two integer \textit{flow variables} $f_{ij}^{ak}$ and $f_{ji}^{ak}$ representing the amount of flow between the cells in the
directions from $i$ to $j$ and from $j$ to $i$, respectively.

The idea of the flow-based model is that a single \textit{source cell} generates a flow equal to the size of the polyomino block, and this flow must be distributed across all cells of the block. If the block consists of several connected components, some cells will be unreachable by the flow, and the corresponding constraints will not be satisfied.

An example of such a flow is shown in Fig.~\ref{Fig:flow_polyomino_example}. 
For the source cell $s$, the outgoing flow equals the size of the polyomino minus one.
For every other cell, the difference between incoming and outgoing flow is equal to one.

\begin{figure}[t]
	\centering
	\begin{tikzpicture}[scale=1.35]
		
		\begin{scope}[every node/.style={rectangle,draw,thick,minimum size=8mm}]
			\node (12) at (2,2) {$s$};
			\node (21) at (1,1) {};
			\node (22) at (2,1) {};
			\node (23) at (3,1) {};
			\node (33) at (3,0) {};
		\end{scope}
		
		\draw [thick,->,>=stealth] (12) edge node[right] {$4$} (22);
		\draw [thick,->,>=stealth] (22) edge node[above] {$1$} (21);
		\draw [thick,->,>=stealth] (22) edge node[above] {$2$} (23);
		\draw [thick,->,>=stealth] (23) edge node[right] {$1$} (33);
		
	\end{tikzpicture}
	\caption{An example of a flow used to verify polyomino connectivity}
	\label {Fig:flow_polyomino_example}
\end{figure}

We do not know in which cell of the grid the source is located, since the exact shape of the block is not known in advance. However, the key observation is that exactly one square of each block lies on the arrow, and we place the source in that
square. Thus, all cells of the polyomino block are divided into \textit{consumers} (those not lying on the arrow) and a single \textit{source} (the one lying on the arrow).

For each block $k$ on each arrow $a$, and for every cell $i \in \mathcal{P}_a$, we introduce a \textit{flow‑supply variable}
\[
F^{ak}_i = 
\begin{cases}
	N_{ak}, & \text{if cell } i \text{ is the source of block } k \text{ on arrow } a,\\
	0, & \text{otherwise}.
\end{cases}
\]

Connectivity of the polyomino blocks is enforced by the following constraints.

\begin{itemize}
	
	\item Total flow generated by the cells on the arrow equals the size of the polyomino block:
	\begin{align}
		\sum_{i \in \mathcal{P}_a} F_{i}^{ak} &= N_{ak} && \forall a \in \mathcal{A}, k = 1,\ldots,K_a.
		\label{constr:source-flow-sum-arrow}
	\end{align}
	We sum over all cells on the arrow because the exact location of the source is not known in advance.
	
	\item Flow may be generated only by a source that belongs to the block:
	\begin{align}
		0 \leq F_{i}^{ak} &\leq M \cdot y_{i}^{ak} && \forall a \in \mathcal{A}, \forall i \in \mathcal{P}_a, k = 1,\ldots,K_a.
		\label{constr:flow-source-cell}
	\end{align}
%	The idea is that if $y^{ak}_i=0$, then $F^{ak}_i=0$.
	
	\item Flow can leave a cell only if it belongs to the block:
	\begin{align}
		f_{ij}^{ak} &\leq M \cdot y_{i}^{ak} &&\forall a \in \mathcal{A}, \forall i \in \mathcal{C}_a, \forall j \in \mathcal{N}(i), k = 1,\ldots,K_a. \label{constr:flow-outbound-cell}
	\end{align}
	
	\item Flow may enter only those cells that belong to the block:
	\begin{align}
		f_{ij}^{ak} &\leq M \cdot y_{j}^{ak} && \forall a \in \mathcal{A}, \forall i \in \mathcal{C}_a, \forall j \in \mathcal{N}(i), k = 1,\ldots,K_a. \label{constr:flow-incoming-cell}
	\end{align}
	
	\item Flow cannot be negative:
	\begin{align}
		f_{ij}^{ak} \;\ge\; 0
		&& \forall a \in \mathcal{A},\ \forall i \in \mathcal{C}_a,\ \forall j \in \mathcal{N}(i),\ k = 1,\ldots, K_a. \label{constr:nonnegative-flow}
	\end{align}
	
		\item Flow conservation for consumer cells $i \notin \mathcal{P}_a$:
	\begin{align}
		\sum_{j \in \mathcal{N}(i)} f_{ji}^{ak}
		- \sum_{j \in \mathcal{N}(i)} f_{ij}^{ak}
		= y_{i}^{ak}
		&& \forall a \in \mathcal{A},\ \forall i \in \mathcal{C}_a \setminus \mathcal{P}_a,\ k = 1,\ldots, K_a. \label{constr:flow-balance-consumers}
	\end{align}
	For a consumer cell, one unit of flow remains in the cell.
	
	\item Flow conservation for a cell $i \in \mathcal{P}_a$ on arrow $a$, which may serve as the source:
	\begin{align}
		\sum_{j \in \mathcal{N}(i) \setminus \mathcal{P}_a} f_{ji}^{ak}
		- \sum_{j \in \mathcal{N}(i) \setminus \mathcal{P}_a} f_{ij}^{ak}
		= y_i^{ak} - F_i^{ak}
		&& \forall a \in \mathcal{A},\ \forall i \in \mathcal{P}_a,\ k = 1,\ldots, K_a. \label{constr:flow-balance-source}
	\end{align}
	For the source cell, the outgoing flow equals $F_i^{ak} - 1$.

\end{itemize}

\subsection{Block evolution variables and constraints}

The final group of constraints ensures that each subsequent block along the arrow adds exactly one square to the previous block, without rotation or reflection. The sizes of the blocks have already been enforced by the basic constraints; what remains is to preserve the shape of the blocks.

For each block $2 \leq k \leq K_a$ on arrow $a \in \mathcal{A}$, and for each admissible two-dimensional translation $(p,q) \in \mathcal{T}_{a}$, we introduce a \textit{translation variable}
\[
t^{ak}_{pq} =
\begin{cases}
	1, & \text{if block } k \text{ on arrow } a \text{ is translated by } p \text{ rows }\\ 
	& \text{and } q \text{ columns relative to block } k-1,\\
	0, & \text{otherwise}.
\end{cases}
\]

An example of a translation is shown in Fig.~\ref{Fig_translation_example}. 
The second block along the arrow contains a copy of the first block (highlighted in blue) that is shifted $2$ rows upward and $3$ columns to the right.

\begin{figure}[t]
	\centering
	\begin{tikzpicture}[scale=0.75]
		
		\foreach \x in {0,...,5}
		\draw (\x+0.5,0.5)--(\x+0.5,3.5);
		
		\foreach \y in {0,...,3}
		\draw (0.5,\y+0.5)--(5.5,\y+0.5);
		
		\draw [semithick,->,>=stealth] (1,1) -- (1,3) -- (4+0.1,3);
		
		\foreach \x in {1}
		\foreach \y in {1}
		\draw (\x-0.25,\y-0.25) -- (\x+0.25,\y-0.25) -- (\x+0.25,\y+0.25) -- (\x-0.25,\y+0.25) -- cycle;
		
		\foreach \x in {2}
		\foreach \y in {1}
		\draw (\x-0.25,\y-0.25) -- (\x+0.25,\y-0.25) -- (\x+0.25,\y+0.25) -- (\x-0.25,\y+0.25) -- cycle;
		
		\foreach \x in {4}
		\foreach \y in {2}
		\draw (\x-0.25,\y-0.25) -- (\x+0.25,\y-0.25) -- (\x+0.25,\y+0.25) -- (\x-0.25,\y+0.25) -- cycle;
		
		\foreach \x in {4}
		\foreach \y in {3}
		\draw [blue] (\x-0.25,\y-0.25) -- (\x+0.25,\y-0.25) -- (\x+0.25,\y+0.25) -- (\x-0.25,\y+0.25) -- cycle;
		
		\foreach \x in {5}
		\foreach \y in {3}
		\draw [blue] (\x-0.25,\y-0.25) -- (\x+0.25,\y-0.25) -- (\x+0.25,\y+0.25) -- (\x-0.25,\y+0.25) -- cycle;
		
	\end{tikzpicture}
	\caption{Example of a translation with $t^{a2}_{-2,3}=1$, i.e., a shift of $2$ rows upward and $3$ columns to the right; the copy of the first block inside the second along the arrow is highlighted in blue}
	\label {Fig_translation_example}
\end{figure}

The following constraints guarantee the correct evolution of the blocks.

\begin{itemize}
	\item Only one translation can be assigned to each block:
	\begin{align}
		\sum_{(p,q) \in \mathcal{T}_{a}} t_{pq}^{ak} &= b_{ak} &&\forall a \in \mathcal{A}, k = 2,\ldots,K_a. \label{constr:unique-translation}
	\end{align}
	
	\item Each square of the block must be translated:
	\begin{align}
		\sum_{(p,q) \in \mathcal{T}^{a}_i} t_{pq}^{ak} &\geq y_{i}^{a,k-1} + b_{ak} - 1 &&\forall a \in \mathcal{A}, \forall i \in \mathcal{C}_a, k = 2,\ldots,K_a,\label{constr:translate-each-square}
	\end{align}
	where $\mathcal{T}^a_i$ is the subset of all admissible translations from cell $i$ in the direction of arrow $a$.
	
	Here, if cell $i$ contains a square of block $k-1$ ($y^{a,k-1}_i = 1$) and block $k$ exists ($b_{ak} = 1$), then that square must have at least one valid translation into block $k$.
	
	\item Correspondence constraints between squares of block $k-1$ and block $k$:
	\begin{align}
		y_{\operatorname{shift}(i,p,q)}^{ak} &\geq y_{i}^{a,k-1} - (1 - t_{pq}^{ak}) &&\forall a \in \mathcal{A}, \forall i \in \mathcal{C}_a, \forall (p,q) \in \mathcal{T}^{a}_i, k = 2,\ldots,K_a. \label{constr:correspondence-squares}
	\end{align}
	
	This constraint preserves the geometric shape of the polyomino blocks across consecutive steps without rotation or reflection.
	If translation $(p,q)$ is chosen for block $k$ ($t^{ak}_{pq} = 1$), the inequality \eqref{constr:correspondence-squares} enforces
	\[ y_{\operatorname{shift}(i,p,q)}^{ak} \geq y_{i}^{a,k-1}, \]
	so a square of block $k-1$ located at cell $i$ must reappear in block $k$ at the translated position $\operatorname{shift}(i,p,q) = i+p\cdot n + q$.
	Otherwise, if $t^{ak}_{pq} = 0$, the constraint reduces to the trivial bound.

\end{itemize}

\subsection{Soundness and completeness of the ILP model}

\begin{theorem}
	Let $\mathcal{E}$ be an instance of the Evolomino puzzle on a board $\mathcal{C}$ with arrow set $\mathcal{A}$, and $S(\mathcal{E})$ be the linear system consisting of constraints
	(\ref{constr:one-block-per-cell})--(\ref{constr:correspondence-squares}),
	with all variables intended to be binary or integral.
	Then an assignment to the cell variables $\{x_i\}_{i\in\mathcal C}$ corresponds to a valid solution of $\mathcal {E}$
	if and only if there exists an integral assignment to the auxiliary variables that satisfies $S(\mathcal{E})$.
\end{theorem}

\begin{proof}
	\textit{$(\Rightarrow)$ Every valid Evolomino solution satisfies the ILP model.}
	
	Given a valid Evolomino solution, set $x_i=1$ on drawn squares and assign values to the auxiliary variables as follows. Index the blocks along each arrow and set the variables $y^{ak}_i$, $b_{ak}$, and $N_{ak}$ accordingly. 
	The constraints \eqref{constr:one-block-per-cell}--\eqref{constr:block-size-increment-leq} follow directly from the definition of blocks and the rules of Evolomino: each square belongs to exactly one block~\eqref{constr:one-block-per-cell}; blocks are ordered along arrows \eqref{constr:block-activation-order}--\eqref{constr:block-order-along-arrow}; adjacent squares belong to the same block \eqref{constr:no-adjacent-blocks-horizontal}--\eqref{constr:no-adjacent-blocks-vertical}; and successive blocks differ in size by one \eqref{constr:block-size-increment-geq}--\eqref{constr:block-size-increment-leq}.
	
	For the flow variables, root a spanning tree of each connected polyomino block at its unique arrow-path cell $s$ and direct all edges away from $s$. Assign to each directed edge $(u,v)$ a flow equal to the number of cells in the subtree rooted at $v$. With this assignment, constraints \eqref{constr:source-flow-sum-arrow}--\eqref{constr:flow-balance-source} are satisfied by standard flow conservation.
	
	For the translation variables, the evolutionary rule guarantees that there exists a unique two-dimensional shift $(p^{*}, q^{*})$ that maps block $k-1$ to a sub-shape of block $k$. By setting $t^{ak}_{p^{*}q^{*}} = 1$ and $t^{ak}_{pq} = 0$ for all $(p,q) \neq (p^{*},q^{*})$, we satisfy constraints \eqref{constr:unique-translation}--\eqref{constr:correspondence-squares}.
	
	\textit{(\(\Leftarrow\)) Every integer-feasible point of the model (\ref{constr:one-block-per-cell})--(\ref{constr:correspondence-squares}) defines a valid solution of the Evolomino puzzle.}
	
	Let $(x, y, b, t, N, f, F)$ be an integer-feasible point. We verify Evolomino rules. 
	Constraint \eqref{constr:one-block-per-cell} ensures that every square $x_i = 1$ belongs to exactly one block, and \eqref{constr:one-square-per-active-block-on-arrow} guarantees that each block contains exactly one square on the pre-drawn arrow. The conditions $b_{a1}=b_{a2}=1$ together with  \eqref{constr:block-activation-order} ensure at least two blocks per arrow.
	Constraints \eqref{constr:no-adjacent-blocks-horizontal}--\eqref{constr:no-adjacent-blocks-vertical} prevent orthogonally adjacent squares from belonging to different blocks, so blocks are geometrically separated.
	
	Connectivity of each polyomino block follows from the flow model: any disconnected component would contain a cell unreachable from the unique source, which contradicts flow conservation \eqref{constr:flow-balance-consumers}--\eqref{constr:flow-balance-source}, since \eqref{constr:flow-source-cell}--\eqref{constr:flow-incoming-cell} confine flow only to the cells of the block ($y^{ak}_{i}=1$).

	Finally, \eqref{constr:unique-translation} selects exactly one translation $(p^{*},q^{*})$, constraints \eqref{constr:translate-each-square}--\eqref{constr:correspondence-squares} together enforce that every square $i$ of block $k-1$ must reappear in block $k$ at the translated position $\operatorname{shift}(i,p^{*},q^{*})$, and \eqref{constr:block-size-increment-geq}--\eqref{constr:block-size-increment-leq} fix $N_{ak}=N_{a,k-1}+1$; hence block \(k\) is exactly the translated copy of block \(k-1\) plus one new square. 
	
	Thus, the integer‑feasible point $(x, y, b, t, N, f, F)$ induces a partition of the drawn squares into disjoint polyomino blocks $B_{ak}=\{i\in\mathcal C : y^{ak}_{i}=1\}$,	
	where each $B_{ak}$ is connected, contains exactly one arrow cell, is obtained from \(B_{a,k-1}\) by a single translation plus one additional square, and not touching other blocks orthogonally. Hence $(x, y, b, t, N, f, F)$ defines a valid Evolomino solution. \qed

\end{proof}

\subsection{Model size}

We estimate how the number of variables and constraints in the ILP model \eqref{constr:one-block-per-cell}--\eqref{constr:correspondence-squares} scales with the size of the Evolomino board, in order to determine the range of instances that can be solved in practice.

For a given rectangular Evolomino board, let $|\mathcal{C}|$ denote the number of cells on the board, and let $K$ be the maximum possible number of polyomino blocks in a solution, where
\[
K = \sum_{a \in \mathcal{A}} K_a.
\]

It is easy to see that $K = O(|\mathcal{C}|)$ and, in fact, $K < |\mathcal{C}|$.
In typical instances, the number of polyomino blocks in a solution is significantly smaller than the board size. 
However, it is technically possible to construct an Evolomino instance whose solution contains a number of blocks that grows linearly with the size of the board (for example, by densely placing short arrows of length 3, each supporting two blocks of sizes 1 and 2, respectively).

We do not attempt to compute the exact number of variables and constraints.
Instead, we estimate their overall magnitude and identify which components contribute the most.

%\subsubsection{Variables.}
\paragraph{Variables.}
Among all variables, the dominant contribution comes from the block variables $y^{ak}_{i}$, the flow variables $f^{ak}_{ij}$, and the translation variables $t^{ak}_{pq}$. 
The size of each of these classes is $O(K \cdot |\mathcal{C}|)$. 
Consequently, the total number of variables in the model is also $O(K \cdot |\mathcal{C}|)$.

%\subsubsection{Constraints.}
\paragraph{Constraints.}
Among constraints, the most significant contributions come from:
\begin{itemize}
	\item the constraints preventing adjacent cells from belonging to different blocks
	\eqref{constr:no-adjacent-blocks-horizontal}--\eqref{constr:no-adjacent-blocks-vertical},
	whose total count is $O(K^2 \cdot |\mathcal{C}|)$;
	
	\item the connectivity constraints for polyomino blocks enforced via flows
	\eqref{constr:source-flow-sum-arrow}--\eqref{constr:flow-balance-source},
	amounting to $O(K \cdot |\mathcal{C}|)$ (with the constant hidden in the big‑$O$ being substantially larger than for other constraints of the same asymptotic order);

	\item the correspondence constraints that maintain square‑by‑square consistency during block evolution along the arrow direction \eqref{constr:correspondence-squares},
	whose total number is $O(K \cdot |\mathcal{C}|^2)$.
\end{itemize}

Thus, the overall number of constraints in the proposed ILP model for the Evolomino puzzle is $O(K \cdot |\mathcal{C}|^2)$, with the most expensive constraints being those that preserve the geometric shape during block evolution.

As an example, the ILP model for the $5 \times 5$ puzzle shown in
Fig.~\ref{Fig_Evolomino_example} contains $509$ variables and $1\,709$
constraints. In contrast, the model for the more complex $10 \times 10$ puzzle from Jonathan Chaffer’s YouTube video~\cite{Pathologic} contains $6\,334$ variables and $55\,057$ constraints.

\section{Puzzle generation}

To construct the dataset for the computational experiments, we developed an algorithm that generates Evolomino puzzles with a unique solution. Its pseudocode is presented in Algorithm~\ref{Alg_Evolomino_generation}.

\begin{algorithm}[p]
	{\footnotesize
		\caption{Evolomino puzzle generation}
		\label{Alg_Evolomino_generation}
		\begin{algorithmic}[1]
			\Procedure{GenerateEvolomino}{$m,n,P$}
			\State $Board \gets$ empty $m\times n$ grid;\ $Tries \gets 0$
			\Comment{We start with an empty grid}
			\While{$\text{Fill}(Board) < P.Target \;\text{and}\; Tries < P.Tries$}
			\If{\Call{TryAddArrow}{$Board, P$} = fail} 
			\State \Call{Backtrack}{$Board$}; $Tries \gets Tries + 1$ \Comment{Undo the last arrow or reset}
			\EndIf
			\EndWhile
			\State \Return \Call{CarveToUnique}{$Board$} \Comment{Final clue reduction and shading}
			\EndProcedure
			
			\Procedure{TryAddArrow}{$Board, P$}
			\State $path \gets [c \gets \text{rand cell}]; \ d \gets \text{rand free dir}$
			\Repeat
			\State $c \gets c + d$; $path.push(c)$ 
			\State $d \gets$ \Call{GetNextDir}{$c, d, Board, P$}  \Comment{Biased random walk}
			\Until{\Call{IsStuck}{$c, Board$} \textbf{or} \Call{Bernoulli}{$P.p\_stop\_arrow$}}
			\If{$|path| < P.min\_arrow$} 
			\State \Return fail 
			\EndIf
			\State $res \gets$ \Call{PlaceBlocks}{$Board, path, \emptyset, 0$}
			\If{$res \neq \emptyset$} 
			\State \Call{Commit}{$Board, res$}; \Return success \Comment{Apply successful placement}
			\EndIf
			\State \Return fail
			\EndProcedure
			
			\Procedure{PlaceBlocks}{$Board, path, blocks, idx$}
			\If{$|blocks| \ge P.min\_blocks$ \textbf{and} \Call{Bernoulli}{$P.p\_stop\_blocks$}} 
			\State \Return $blocks$ 
			\EndIf
			\State $Cands \gets \Call{GetValidIndices}{path, idx, blocks}$
			\For{\textbf{each} $j \in \Call{WeightedShuffle}{Cands}$} \Comment{Weighted random order}
			\State $Shapes \gets$ \Call{GetEvolutions}{$blocks.last$}
			\For{\textbf{each} $S \in$ \Call{Shuffle}{$Shapes$}}
			\For{\textbf{each} $anchor \in$ \Call{Shuffle}{$S$}}	
			\State $B \gets \Call{Align}{S, anchor, path[j]}$ \Comment{Align $anchor$ to $path[j]$}
			\If{\Call{CanPlace}{$B, Board$}}
			\State $res \gets \Call{PlaceBlocks}{Board, path, blocks \cup \{B\}, j}$
			\If{$res \neq \emptyset$}
			\State \Return $res$ 
			\EndIf
			\State \Call{Undo}{$B, Board$} \Comment{Backtracking: revert failed placement}
			\EndIf
			\EndFor
			\EndFor
			\EndFor
			\State \Return $\emptyset$
			\EndProcedure
			
			\Procedure{CarveToUnique}{$Board$} \Comment{Greedy clue removal and shading}
			\State $Sol \gets$ \Call{ShadeAllFreeCells}{$Board$} \Comment{Fix a unique solution}
			\State $Cells \gets \Call{Shuffle}{\text{AllCells}}$ \Comment{Random order to avoid bias}
			\ForAll{$c \in Cells$}
			\State $Old \gets Board[c]$
			\If{$Old \in \{\text{SHADED}, \text{SQUARE}\}$}
			\State $Board[c] \gets \text{EMPTY and UNSHADED}$
			\If{\Call{HasDifferentSolution}{$Board, Sol$}}
			\State $Board[c] \gets Old$ \Comment{Revert if an alternative solution exists}
			\EndIf
			\EndIf
			\EndFor
			\State \Return $Board$ \Comment{Resulting puzzle with a unique solution}
			\EndProcedure
		\end{algorithmic}
	}
\end{algorithm}

The algorithm starts from an empty rectangular board and gradually fills it while attempting to reach a prescribed density (the total length of arrows and the total size of polyomino blocks). At each step, a new arrow is drawn using a biased random walk, where continuing straight is preferred over turning. The process stops when the arrow cannot grow further or when a probabilistic stopping criterion is triggered, implemented throughout the algorithm via repeated Bernoulli trials. Failed arrows are rolled back, and a new attempt is made.

For every admissible arrow, the block‑placement procedure is executed.
It first selects a random anchor position on the arrow for the initial block, prioritizing positions closer to the beginning of the arrow (weighted random sampling by Efraimidis and Spirakis~\cite{Efraimidis2006}).
From this anchor, a random polyomino block is grown.
Subsequent blocks are placed afterward, with their number treated as a random variable, again drawn from a geometric distribution.
For each candidate, the procedure enumerates all feasible evolution shapes, anchor positions on the arrow, and anchor‑cell choices within the polyomino.
Conflicts trigger a rollback and a new attempt.
A successful sequence of blocks is committed to the board.

In the final stage, after all arrows and blocks have been placed and the required density has been achieved, the algorithm shades all remaining empty cells, thereby fixing a unique solution.
It then iterates through the cells in random order (classic Fisher–Yates shuffle~\cite{Knuth1997}), greedily attempting to remove clue squares and shaded cells.
If a removal does not introduce an alternative solution, as verified by calling the ILP solver, where the current solution $x^\ast$ is excluded by the additional constraint
\[\sum_{i \in \mathcal{C}} |x_i - x_i^\ast| \ge 1,\]
the change is accepted; otherwise, it is rolled back.
The resulting puzzle contains a minimal set of clues that still uniquely determines the constructed solution.

\section{Computational experiments}

Using the generator described in Algorithm~\ref{Alg_Evolomino_generation}, we constructed a custom benchmark dataset (available at~\cite{EvolominoDataset2026}) that includes:
\begin{itemize}
	\item $500$ puzzles with unique solutions, ranging in size from $5 \times 5$ to $14 \times 14$, with $50$ instances for each grid size;
	\item $200$ puzzles with guaranteed, but not necessarily unique, solutions, ranging from $15 \times 15$ to $18 \times 18$, again evenly distributed across grid sizes.
\end{itemize}

For larger instances, we had to replace the $\textsc{CarveToUnique}$ procedure, which guarantees solution uniqueness, with random clue removal and obstacle placement, since repeated calls to the ILP solver became computationally expensive.

The proposed ILP model was evaluated using the open-source software suite Google OR-Tools v9.15~\cite{OR-Tools}. All experiments were performed on a machine equipped with an Intel Core i7--13620H CPU and 16\,GB of DDR5--5400 RAM, running Windows~11~Pro.

We conducted two series of experiments. In the first, using $50$ puzzle instances of size $10 \times 10$, we compared four solvers available out-of-the-box in OR-Tools:
\begin{itemize}
	\item CP-SAT: an award-winning constraint programming and SAT solver~\cite{cpsatlp};
	
	\item SCIP: a high-performance academic mixed-integer programming solver~\cite{SCIP10};
	
	\item CBC: an open-source branch-and-cut solver from the COIN-OR initiative~\cite{Lougee-Heimer2004CBC};
	
	\item BOP: a specialized SAT-based solver for pure Boolean problems~\cite{OR-Tools}.
\end{itemize}

This comparison was performed to identify the most efficient solver for our ILP model of the Evolomino puzzle.

The performance profile of this experiment is presented in Fig.~\ref{Fig_performance_profile}, with each point representing two instances. CP-SAT emerges as the clear winner, being the only solver to successfully handle all 50 benchmark instances within the 180-second time limit. In comparison, SCIP solved 49 instances, while both BOP and CBC completed 48.

Within a one-second cutoff, CP-SAT solved 43 out of 50 instances, BOP solved 22, SCIP 14, and CBC only 9. A notable outlier for CP-SAT is \textit{sample39}, which required 12 seconds to solve. This instance proved particularly challenging: under the 3-minute limit, only CP-SAT and BOP managed to solve it at all.

\begin{figure}[t]
	\centering
	\begin{tikzpicture}
		\begin{axis}[
			width=\textwidth,
			height=9.5cm,
			xmode=log,
			log basis x={10},
			xlabel={Computation time (s)}, 
			ylabel={Solved instances (\%)},
			xticklabels={0.01, 0.1, 1, 10, 100},
			ymin=-0.05, ymax=1.05,
			ytick={0.0,0.1,0.2,0.3,0.4,0.5,0.6,0.7,0.8,0.9,1.0},
			yticklabel={\pgfmathparse{\tick*100}\pgfmathprintnumber{\pgfmathresult}\%},
			grid=major,
			grid style={dashed,gray!40},
			legend pos=south east,
			thick,
			mark options={scale=1.25}
			]
			
			%CP-SAT
			\addplot[color=blue, mark=*] coordinates {
				(142, 0.04) (174, 0.08) (220, 0.12) (271, 0.16) (326, 0.20) (391, 0.24) (445, 0.28) (500, 0.32) (523, 0.36) (546, 0.40) (554, 0.44) (585, 0.48) (599, 0.52) (620, 0.56) (750, 0.60) (771, 0.64) (786, 0.68) (797, 0.72) (842, 0.76) (856, 0.80) (952, 0.84) (1075, 0.88) (1110, 0.92) (1172, 0.96) (12772, 1.00)
			};
			\addlegendentry{CP-SAT}

			% BOP
			\addplot[color=orange, mark=halfcircle*] coordinates {
				(160, 0.04) (191, 0.08) (273, 0.12) (511, 0.16) (531, 0.20) (608, 0.24) (733, 0.28) (826, 0.32) (862, 0.36) (925, 0.40) (977, 0.44) (1079, 0.48) (1151, 0.52) (1413, 0.56) (1492, 0.60) (1797, 0.64) (1968, 0.68) (2190, 0.72) (2321, 0.76) (2538, 0.80) (3608, 0.84) (6311, 0.88) (15432, 0.92) (78357, 0.96)
			};
			\addlegendentry{BOP}
			
			% SCIP
			\addplot[color=teal, mark=pentagon*] coordinates {
				(142, 0.04) (271, 0.08) (328, 0.12) (462, 0.16) (608, 0.20) (865, 0.24) (954, 0.28) (1013, 0.32) (1163, 0.36) (1288, 0.40) (1609, 0.44) (2016, 0.48) (2125, 0.52) (2307, 0.56) (2399, 0.60) (2595, 0.64) (2741, 0.68) (3169, 0.72) (3616, 0.76) (4612, 0.80) (6295, 0.84) (6518, 0.88) (7402, 0.92) (7994, 0.96) (10393, 0.98)
			};
			\addlegendentry{SCIP}
			
			% CBC
			\addplot[color=purple, mark=triangle*] coordinates {
				(570, 0.04) (672, 0.08) (723, 0.12) (898, 0.16) (1160, 0.20) (1769, 0.24) (1825, 0.28) (1919, 0.32) (3248, 0.36) (3417, 0.40) (3738, 0.44) (4074, 0.48) (4276, 0.52) (4575, 0.56) (5419, 0.60) (6059, 0.64) (6722, 0.68) (7150, 0.72) (10034, 0.76) (16748, 0.80) (31430, 0.84) (39189, 0.88) (44553, 0.92) (108757, 0.96)
			};
			\addlegendentry{CBC}

		\end{axis}
	\end{tikzpicture}
	\caption{Performance profile of the solvers on a $10 \times 10$ grid (180-second time limit)}
	\label{Fig_performance_profile}
\end{figure}

\begin{figure}[p]
	\centering
	\small
	\begin{tikzpicture}
		\begin{axis}[
			width=\textwidth,
			height=9.5cm,
			xlabel={Grid size $n \times n$},
			ylabel={Computation time (s)},
			ymode=log,
			log basis y={10},
			ymin=10,
%			yminorticks=false,
			yticklabels={0.01, 0.1, 1, 10, 100},
			%			xmode=log,
			%			log basis x={2},
			xmin=4.5, xmax=18.5,
			xtick={5, 6, 7, 8, 9, 10, 11, 12, 13, 14, 15, 16, 17, 18},
			legend pos=south east,
			grid=major,
			grid style={dashed,gray!40},
			enlarge x limits=false,
			]

			\addplot[name path=q3, draw=none, forget plot] coordinates {
				(5, 32) (6, 61) (7, 111) (8, 234) (9, 438) 
				(10, 842) (11, 1736) (12, 3301) (13, 5327) (14, 9741) (15, 22832) (16, 37457) (17, 74549) (18, 80257)
			};
			
			\addplot[name path=q1, draw=none, forget plot] coordinates {
				(5, 25)  (6, 41) (7, 79) (8, 127) (9, 223) 
				(10, 415) (11, 626) (12, 1681) (13, 2127) (14, 4727) (15, 9662) (16, 10409) (17, 20165) (18,31961)
			};
			
			\addplot[fill=blue!20, forget plot] fill between[of=q3 and q1];
			
			\addplot[color=blue, thick, mark=*, mark size=1.5pt] coordinates {
				(5, 27) (6, 53) (7, 95) (8, 175.5) (9, 334) 
				(10, 597.5) (11, 1246) (12, 2547) (13, 3761) (14, 7643.5) (15, 14952.5) (16, 23543.5) (17, 39130.5) (18,54739.5)
			};
			\addlegendentry{CP-SAT}
			
			\addplot[color=blue,only marks, mark size=1pt]
			coordinates {
				(5, 45) (5, 47) (5, 57) (6,175) (7,186) (7,238) (8,444) 
				(10, 12772) (12, 7681) (12, 17892) (15, 45834) (15, 46795) (15, 52699) (15, 53127) (15, 63747) (15, 68324) (15, 104406) (16, 89658) (16, 91876) (16, 116200) (17, 158943) (17, 170273) (17, 173732) (17, 185554) (17, 271873) (18, 187494) (18, 310658)
			};

		\end{axis}
	\end{tikzpicture}
	\caption{CP-SAT performance: median (line), IQR (shaded), and outliers (dots)}
	\label{Fig_median_runtime_CP_SAT}
\end{figure}

\begin{table}[p]
	\centering
	\caption{ILP model complexity and CP-SAT solver runtime on the benchmark dataset~\cite{EvolominoDataset2026} (50 instances per grid size)}
	\label{Tab_computational_experiments}
	\begin{tabular}{|r|r|r|r|r|r|}
		\hline
		\multirow{2}{*}{~Grid size~} & \multicolumn{2}{|c|}{Problem size (Med.)} & \multicolumn{3}{|c|}{Running time (s)} \\ \cline{2-6}
		& ~Variables~ & ~Constraints~ &  ~~~~~$Q1$~~~~~ & ~~Median~~ & ~~~~~$Q3$~~~~~\\
		\hline
		$5 \times 5$~ & $450$~ & $1\,429$~ & $0.025$~ & $0.027$~ & $0.032$~ \\ 
		$6 \times 6$~ & $876$~ & $3\,615$~ &  $0.041$~ & $0.053$~ & $0.061$~ \\ 
		$7 \times 7$~ & $1\,525$~ & $7\,556$~ &  $0.079$~ & $0.095$~ & $0.111$~ \\ 
		$8 \times 8$~ & $2\,471$~ & $14\,326$~ &  $0.127$~ & $0.175$~ & $0.234$~ \\ 
		$9 \times 9$~ & $3\,819$~ & $26\,729$~ &  $0.223$~ & $0.334$~ & $0.438$~ \\ 
		$10 \times 10$~ & $5\,352$~ & $41\,905$~ &  $0.415$~ & $0.597$~ & $0.842$~ \\ 
		$11 \times 11$~ & $8\,080$~ & $74\,248$~ &  $0.626$~ & $1.246$~ & $1.736$~ \\ 
		$12 \times 12$~ & $11\,287$~ & $129\,946$~ &  $1.681$~ & $2.547$~ & $3.301$~ \\ 
		$13 \times 13$~ & $14\,940$~ & $192\,501$~ &  $2.127$~ & $3.761$~ & $5.327$~ \\ 
		$14 \times 14$~ & $20\,351$~ & $289\,132$~ &  $4.727$~ & $7.643$~ & $9.741$~ \\
		$15 \times 15$~ & $27\,648$~ & $459\,073$~ &  $9.662$~ & $14.952$~ & $22.832$~ \\
		$16 \times 16$~ & $31\,893$~ & $525\,008$~ &  $10.409$~ & $23.543$~ & $37.457$~ \\
		$17 \times 17$~ & $39\,675$~ & $740\,220$~ &  $20.165$~ & $39.130$~ & $74.549$~ \\
		$18 \times 18$~ & $53\,307$~ & $1\,121\,290$~ &  $31.961$~ & $54.739$~ & $80.257$~ \\
		\hline
	\end{tabular}
\end{table}

In the second series of tests, the CP-SAT solver that demonstrated the best performance was evaluated on the entire dataset of $700$ benchmarks (comprising $500$ instances with unique solutions and $200$ with guaranteed solutions). The results of these computational experiments are shown in Fig.~\ref{Fig_median_runtime_CP_SAT} and in more detail in Table~\ref{Tab_computational_experiments}. The solid line indicates the median runtime, the shaded area represents the interquartile range (IQR), and individual points denote Tukey outliers ($Q_1 - 1.5\,\mathrm{IQR}$ and $Q_3 + 1.5\,\mathrm{IQR}$~\cite{Dekking2005}).

The plot demonstrates that solving time grows exponentially with the size of the puzzle grid and, consequently, with the number of linear constraints in the ILP model. For smaller grids up to $10 \times 10$, the median runtime remains below one second; for grids up to $13 \times 13$, it stays within five seconds. Starting from $15 \times 15$, the IQR expands and outliers become more frequent, indicating that the instances become less homogeneous and the solver's behavior less predictable. We stopped at puzzles of size $18 \times 18$, where the median runtime per instance ($54$ seconds) still remained under one minute.

Given the complexity of the ILP model and the cubic growth $O(K \cdot |\mathcal{C}|^2)$ of the number of linear constraints relative to the grid size, these results are quite promising and highlight the efficiency of modern state-of-the-art solvers. In particular, for successfully solved $18 \times 18$ puzzles, the model contained approximately $50\,000$ variables and about one million linear constraints.

\section{Conclusion}

Evolomino is a relatively new pencil-and-paper logic puzzle that remains largely unexplored from a theoretical computer science perspective. To the best of our knowledge, the only existing result is a proof establishing that determining whether a given Evolomino instance has at least one valid solution is NP-complete and that counting the number of distinct solutions is $\texttt{\#}$P-complete~\cite{Nikolaev2025}.

In this paper, we propose a computational approach to solving the puzzle within the transform-and-conquer paradigm by formulating the Evolomino rules as a set of constraints for an integer linear programming model. Furthermore, we introduce an algorithm for generating random Evolomino instances, where solution uniqueness is guaranteed via the proposed ILP framework.

Computational experiments conducted on a custom benchmark dataset demonstrate that the state-of-the-art CP-SAT solver from the OR-Tools library successfully solves puzzle instances of size up to $10 \times 10$ within one second, and up to $18 \times 18$ within one minute.

Future research could explore alternative approaches to solving the Evolomino puzzle, such as backtracking, constraint programming, SAT-based models~\cite{Weber2005}, or a reduction to the exact cover problem and applying the DLX (Dancing Links) algorithm~\cite{Knuth2000}, which is particularly effective for placing polyominoes on a grid.

\begin{credits}
\subsubsection{\ackname} The research of Andrei V. Nikolaev was supported by Alfa‑Bank JSC under the ``Alfa‑Future: Grants for Faculty'' program. 
We thank the anonymous reviewers for their insightful and constructive comments, which helped improve the clarity and presentation of this paper.

\subsubsection{\discintname}
The authors have no competing interests to declare that are
relevant to the content of this article

%It is now necessary to declare any competing interests or to specifically
%state that the authors have no competing interests. Please place the
%statement with a bold run-in heading in small font size beneath the
%(optional) acknowledgments\footnote{If EquinOCS, our proceedings submission
%system, is used, then the disclaimer can be provided directly in the system.},
%for example: The authors have no competing interests to declare that are
%relevant to the content of this article. Or: Author A has received research
%grants from Company W. Author B has received a speaker honorarium from
%Company X and owns stock in Company Y. Author C is a member of committee Z.
\end{credits}

\bibliographystyle{splncs04}
\bibliography{Nikolaev_Myasnikov_arXiv_2026}

\end{document}